\documentclass[11pt]{article}
\usepackage{latexsym}
\usepackage{amsfonts}
\usepackage{graphicx}
\usepackage{epstopdf}
\usepackage{xcolor}
\usepackage[sort]{natbib}
\usepackage{lineno}
\usepackage{rotating}

\newtheorem{theorem}{Theorem}
\newtheorem{lemma}{Lemma}
\newtheorem{corol}{Corollary}
\newtheorem{property}{Property}

\newcommand{\be}{\begin{equation}}
\newcommand{\ee}{\end{equation}}

\begin{document}
\title{\bf Multi-Facility Location Models Incorporating Multipurpose Shopping Trips}
\author{Pawel Kalczynski\\
	College of Business and Economics\\
	California State University-Fullerton\\
	Fullerton, CA 92834.\\e-mail: pkalczynski@fullerton.edu
	\and	Zvi Drezner\footnote{Corresponding author. e-mail:zdrezner@fullerton.edu.} \\
	College of Business and Economics\\
	California State University-Fullerton\\
	Fullerton, CA 92834.\\e-mail: zdrezner@fullerton.edu	
	\and
	Morton O'Kelly\\
	Department of Geography\\
	Ohio State University\\
	Columbus, OH 43210.\\
	email: okelly.1@osu.edu
}

\date{}
\maketitle


\bigskip	
\begin{abstract}
This paper continues to develop and explore the impact of multipurpose trips on retail location. We develop the model of locating multiple competing facilities of a chain where several competing facilities exist in the area. There may be some existing facilities of the same chain as well.
The addition of  multiple new outlets can cause cannibalization of existing sales, but this effect is mitigated by selecting good locations, and the total market share captured by the chain increases. The introduction of multipurpose trips enhances the total market share of the location decision maker. Since in reality many customers combine a visit to more than one facility in one shopping trip, the model predicts the expected market share captured more accurately. Therefore, the selected locations for new facilities are more accurate as well.	
\end{abstract}

	\renewcommand{\baselinestretch}{2}
	\renewcommand{\arraystretch}{0.5}
	\large
	\normalsize

\section{Introduction}

Competitive facilities location problems attempt to find the locations for one or more new facilities among existing competing facilities that maximize the captured market share. The facilities attract demand generated by customers in the area. In most applications profit is increasing when market share captured increases. Therefore, the common objective is to maximize the market share captured by the new facilities. Multipurpose (MP) shopping models consider the possibility that some customers incorporate in one trip two (or more) competing facilities.

There are two main approaches to model multipurpose shopping trips. In one approach \citep{O81,O83a,O83b,LSS04,M83,MG86}, investigated in this paper, a customer plans a trip to buy a desired product and combines in one trip purchasing a second product in a different facility, which is termed a ``cluster". For example, purchasing coffee in a coffee shop and patronizing a grocery store on the same trip. This may influence the customer's choice of the coffee shop to visit. Customers visit two facilities with a given proportion $\pi$, and  a proportion $1-\pi$ of customers visit only one facility for the desired product. For a recent review of multipurpose shopping see \citet{ZOK23,DOD22,O09}.
A different model \citep{MEL19,MMF23}, also termed multipurpose shopping trip, is to visit several facilities that sell the desired product, and go back to the selected facility to purchase it.

The bulk of competitive facilities models consider the location of facilities assuming a single purpose trip without the possibility of multipurpose trips. There are many approaches for estimating the market share captured by a facility located among facilities that offer the same product.  For full details see \citet{TammyChap,DE23} and the references therein. The most commonly used model is the gravity model proposed by \citet{R31} which was applied in many papers  \citep[for example,][]{NC74,Bell98,JM79,Dr94b,ABK1,ABK2,Huff64,Huff66,FH13,SHFP09}.
 There are other models for estimating the market share, \citep[for example, ][]{LT84,DD96,DDK11,Dr82,ReV86}, which are beyond the scope of this paper and proposed as future research.

Analysts of retail trade areas have used two different approaches to measurement: one is to query a spatially stratified random sample of residents about their typical shopping behavior. This requires a lot of effort. A more efficient approach is to conduct a choice-based sampling of people as they patronize a particular store, which of course gets a large sample from that store, but requires work to make inferences about the population at large \citep[see][]{O99,DO08}. Interestingly, the first approach is suited to multiple purpose shopping because we can frame the survey to include packages of intended trips. The choice-based approach on the other hand is hampered by the fact that the shopper may not yet have completed any other stops or may not yet know of their intention to proceed to further shopping. This motivates the need for a theoretical model that captures the intended or potential array of ancillary locations that are available to the shopper during their trip. We can tell from casual observation that such chaining must occur because of the prevalence of multiple similar clusters of chains that are designed to capture multipurpose trip makers. For example, coffee shops and grocery stores, cell phone providers, and other convenience oriented services.

In the gravity model we assume that there are $n$ demand points in the area, each with a buying power for the particular product of $b_i$ for $1\le i\le n$, and $p$ competing facilities with attractiveness $A_j$ for $1\le j\le p$. The distance between demand point $i$ and facility $j$ is $d_{ij}$. We assume that the proportion of customers that patronize a certain facility is proportional to the facility's attractiveness and declines as the distance increases by a distance decay function $f(d)$. The expected market share $M_j$ captured by facility $j$ is:
\begin{equation}
M_j=\sum\limits_{i=1}^nb_i\frac{A_jf(d_{ij})}{\sum\limits_{k=1}^pA_kf(d_{ik})}
\end{equation}
Note that $\sum\limits_{j=1}^pM_j=\sum\limits_{i=1}^nb_i$, meaning that the competing facilities attract all the available buying power (the sum of the proportions for each demand point is 1).

Commonly investigated distance decay functions (for some parameter $\lambda$) are power decay,  $f(d)=\frac{1}{d^\lambda}$ \citep{Huff64,Huff66}, and exponential decay, $f(d)=e^{-\lambda d}$ \citep{W76a,H81a}. The parameter $\lambda$ depends on the retail category. It may be different for grocery stores than for shopping malls. Note that for exponential decay the parameter $\lambda$ also depends on the units of the distance measure. If $\lambda$ is applied for distances in miles, then $\frac{\lambda}{1.609}$ should be applied if distances are measured in kilometers. \citet{Dr06a} found that the exponential decay function fits real data for shopping malls better than power decay.

The paper is organized as follows. In Section \ref{sec1} we describe the proposed model, and its formulation is detailed in Section \ref{sec2}. In Section \ref{sec3} computational experiments are reported on 72 instances applying the power decay function, and the same 72 instances applying the exponential decay for a total of 144 experiments. We summarize the paper and propose ideas for future research in Section~\ref{sec4}.

\section{\label{sec1}The Proposed Model}
Most papers investigating multipurpose trips consider the location of one competing facility for the desired product. In this paper we investigate locating a chain of facilities that offer the same desired product. There are $\hat p$ competing facilities in the area, and $p\ge 1$ new chain facilities need to be located. The objective is to maximize the market share captured by the chain of new facilities.

It is possible that $0\le q\le \hat p-1$ of the $\hat p$ existing facilities belong to the same chain. The formulation of the market share captured by the chain is not changed by much. We add the $q$ existing facilities to create a chain of $p+q$ ``new" facilities, and retain only $\hat p-q$ existing facilities. However, in the optimization phase only the original $p$ new facilities are variables and the transferred $q$ facilities are fixed. Therefore, regardless of the value of $q$ there are $2p$ variables (the $x$ and $y$ of each location of the new facility) in the optimization phase, and the optimization phase is expected to perform with the same effectiveness for any $q\ge 0$. In the computational experiments we assume for simplicity that $q=0$ because the expected effort for solving the problem does not depend on $q$. Also, testing several values of $q$ will significantly increase the number of test instances. We now have 144 test instances which is definitely sufficient to assess the model.

\section{\label{sec2}Estimating the Captured Market Share}

The estimated market share is calculated by the gravity model \citep{R31,Huff64,Huff66}. Other models will require formulas for the market shares captured by a single purpose trip and those for a multipurpose trip to replace these expressions in Section~\ref{formul}.
 Investigating other models for estimating market share is beyond the scope of this paper and is suggested for future research.

\subsection{Definitions}

\begin{tabular}{cl}
	$n$&The number of demand points.\\
	$b_i$&The buying power at demand point $i$ dedicated to the desired product.\\
	$d_i(X)$&The distance between demand point $i$ and location $X$.\\
	$d(X,Y)$&The distance between locations $X$ and $Y$.\\
	$\hat p$&The number of existing competing facilities for the desired product.\\
	$\hat X_j$&The location of competing facility $j$ for $j=1,\ldots, \hat p$.\\
	$A_j$&Attractiveness of competing facility $j$ for $j=1,\ldots, \hat p$.\\ 
	$p'$&The number of clusters of facilities who sell a different product.\\
	$Y_m$& The location of cluster $m$ for $m=1,\ldots,p'$.\\
	$A'_m$&Attractiveness of cluster $m$ for $m=1,\ldots, p'$.\\
	$p$&The number of new chain facilities to be located.\\
	$X_j$&The unknown location of new facility $j$ for $j=1,\ldots,p$.\\
	$X$&$\{X_j, j=1,\ldots,p\}$.\\
	$\hat A_j$&The attractiveness of new facility $j$.\\ 
	$f(d)$&The distance decay function by the travel distance.\\
	$\pi$&The proportion of multipurpose trips.\\
\end{tabular}	
\medskip

\subsection{\label{formul}Formulation}

We assume symmetric distances, i.e., $d(Y,X)=d(X,Y)$. Therefore, the order of the visits to the two facilities results in the same total distance. The formulations  for locating $p$ facilities are extensions of the formulations in \citet{ZOK23,DOD22} for locating one facility.
The market share captured by the chain of $p$ facilities is defined as $M(X)$. The proportion of the market share captured, reported in the computational experiments, is $\frac{M(X)}{\sum\limits_{i=1}^n b_i}$.

For Power decay the parameter $\lambda=2$ is applied. A distance correction  $\alpha=\frac{24}{\sum\limits_{i=1}^nb_i}$ \citep{DD97} is included. The derivation of $\alpha$ is detailed in \citet{ZOK23}. The captured market share is:

\begin{eqnarray}
M(X)&=&\pi  \sum\limits _{i=1}^n \frac{b_i \sum\limits _{m=1}^{p'} A'_m \sum\limits _{j=1}^p \frac{\hat A_j}{\left(d(X_j,Y_m)+\sqrt{
			\alpha b_i+d_i^2(X_j)}+\sqrt{\alpha b_i+d_i^2(Y_m)}\right)^2}}{\sum\limits _{m=1}^{p'} A'_m \sum\limits _{j=1}^p \frac{\hat A_j}{\left(d(X_j,Y_m)+\sqrt{\alpha b_i+d_i^2(X_j)}+\sqrt{\alpha b_i+d_i^2(Y_m)}\right)^2}~+~C_i^{(2)}}\nonumber\\&+&(1-\pi ) \sum\limits _{i=1}^n \frac{b_i \sum\limits _{j=1}^p \frac{\hat A_j}{ \alpha b_i+d_i^2(X_j)}}{\sum\limits _{j=1}^p \frac{\hat A_j}{\alpha b_i+d_i^2(X_j)}~+~C_i^{(1)}}
\end{eqnarray}
where
$$C_i^{{(1)}}=\sum _{k=1}^{\hat p} \frac{A_k}{d_i^2(\hat X_k)+\alpha b_i}$$$$C_i^{{(2)}}=\sum _{k=1}^{\hat p} \sum _{m=1}^{p'} \frac{A_k A_m}{\left(\sqrt{d^2(\hat X_k,Y_i)+\alpha b_i}+d(\hat X_k,Y_m)+\sqrt{\alpha b_i+d_i^2(Y_m)}\right)^2}~.$$

For Exponential decay and a given $\lambda$, the captured market share is:

\begin{eqnarray}
M(X)&=&\pi  \sum\limits _{i=1}^n \frac{b_i \sum\limits _{m=1}^{p'} A'_m \sum\limits _{j=1}^p \hat A_je^{-\lambda  \left(d\left(X_j,Y_m\right)+d_i\left(X_j\right)+d_i\left(Y_m\right)\right)}}{\sum\limits _{m=1}^{p'} A'_m \sum\limits _{j=1}^p \hat A_je^{-\lambda  \left(d\left(X_j,Y_m\right)+d_i\left(X_j\right)+d_i\left(Y_m\right)\right)}~+~C_i^{(2)}}\nonumber\\&+&(1-\pi ) \sum\limits _{i=1}^n \frac{b_i \sum\limits _{j=1}^p \hat A_je^{-2 \lambda  d_i\left(X_j\right)}}{\sum\limits _{j=1}^p \hat A_je^{-2 \lambda  d_i\left(X_j\right)}~+~C_i^{(1)}}
\end{eqnarray}
where
$$C_i^{{(1)}}=\sum _{k=1}^{\hat p} A_k e^{-2 \lambda  d_i\left(\hat X_k\right)}$$
$$C_i^{{(2)}}=\sum\limits_{k=1}^{\hat p} \sum\limits_{m=1}^{p'} A_k A_m e^{-\lambda  \left(d_i(\hat X_k)+d\left(\hat X_k,Y_m\right)+d_i\left(Y_m\right)\right)}~.$$

\subsection{Solution Method}
\citet{DOD22} solved the $p=1$ instances optimally, within a relative accuracy of  $\epsilon=10^{-5}$, by the BSSS algorithm \citep{HPT81}. The largest problem of $n=20,000$ demand points was solved within 70 minutes. This can be done only for the $p=1$ instances because for $p>1$ there are more than two variables. Problems with $p=2$ and 3 may be solved with the BCSC algorithm \citep{SS10}, but such  large problems are not expected to be solved in a manageable run time.
	We therefore solved the problems with up to $p=20$ new facilities by a local non-linear solver: Sparse Nonlinear Optimizer \citep[SNOPT,][]{GMS05}  without guaranteeing optimality. Each problem was solved from 100 random starting solutions and the best result reported.
	All the results for $p=1$ were identical to those reported in \citet{DOD22} and thus are optimal. Note that for $p=1$ an optimal solution can be found by BSSS, even if some chain facilities already exist in the area, because there are only two variables in the formulation. This problem is also suggested for future research.

\begin{table}[ht!]
		\caption{\label{n100}Results for $n=100$}
\begin{center}
\begin{tabular}{|c||c|c|c|c|c|c|}						
	\hline						
	$p$&$\pi=0.0$&$\pi=0.2$&$\pi=0.4$&$\pi=0.6$&$\pi=0.8$&$\pi=1.0$\\						
	\hline						
	\multicolumn{7}{|c|}{Proportion of Market Share Captured (Power Decay)}\\						
	\hline						
	1&	0.12407&	0.11856&	0.12073&	0.12874&	0.13734&	0.14739\\
	2&	0.22052&	0.21495&	0.21736&	0.22653&	0.23767&	0.25286\\
	3&	0.30230&	0.29850&	0.29748&	0.30469&	0.31851&	0.33431\\
	4&	0.37155&	0.36315&	0.35791&	0.36412&	0.38171&	0.40002\\
	5&	0.42093&	0.41405&	0.41002&	0.41388&	0.43309&	0.45335\\
	10&	0.59494&	0.58386&	0.57673&	0.57810&	0.59549&	0.62019\\
	15&	0.69158&	0.68087&	0.67265&	0.66956&	0.68385&	0.70839\\
	20&	0.75311&	0.74198&	0.73320&	0.72986&	0.73943&	0.76326\\
	\hline\multicolumn{7}{|c|}{Run Time in Minutes for all 100 Runs (Power Decay)}\\\hline						
	1&	0.08&	0.09&	0.20&	0.24&	0.28&	0.27\\
	2&	0.17&	0.17&	0.33&	0.85&	0.78&	0.86\\
	3&	0.31&	0.32&	0.43&	1.10&	1.22&	1.11\\
	4&	0.48&	0.48&	0.88&	1.79&	1.85&	1.54\\
	5&	0.61&	0.69&	1.14&	1.92&	2.14&	1.96\\
	10&	1.62&	1.71&	2.47&	4.23&	5.94&	4.98\\
	15&	2.64&	2.74&	3.89&	7.27&	10.37&	9.11\\
	20&	3.83&	3.93&	6.40&	10.61&	16.09&	15.00\\
	\hline						
	\multicolumn{7}{|c|}{Proportion of Market Share Captured (Exponential Decay)}\\						
	\hline						
	1&	0.14471&	0.13564&	0.12804&	0.15121&	0.17893&	0.20665\\
	2&	0.24657&	0.24338&	0.24376&	0.27880&	0.31616&	0.35360\\
	3&	0.33396&	0.33774&	0.35617&	0.38301&	0.42683&	0.47588\\
	4&	0.41590&	0.41760&	0.42758&	0.46184&	0.51326&	0.57338\\
	5&	0.46944&	0.46846&	0.48230&	0.51094&	0.56957&	0.63149\\
	10&	0.65250&	0.65485&	0.66533&	0.68614&	0.71976&	0.77913\\
	15&	0.75590&	0.74913&	0.75109&	0.76352&	0.78986&	0.83811\\
	20&	0.81296&	0.80714&	0.80532&	0.81195&	0.83145&	0.87258\\
	\hline						
	\multicolumn{7}{|c|}{Run Time in Minutes for all 100 Runs (Exponential Decay)}\\						
	\hline						
	1&	0.18&	0.13&	0.38&	0.46&	0.42&	0.50\\
	2&	0.43&	0.25&	0.67&	0.86&	0.90&	0.86\\
	3&	0.77&	0.49&	1.09&	1.30&	1.32&	1.28\\
	4&	1.30&	0.83&	1.48&	1.54&	1.80&	1.73\\
	5&	0.80&	1.00&	1.60&	2.20&	2.35&	2.17\\
	10&	2.78&	3.28&	3.90&	5.49&	6.91&	6.69\\
	15&	6.18&	7.76&	8.94&	10.94&	13.68&	14.09\\
	20&	14.27&	15.61&	15.87&	19.58&	22.80&	24.37\\
\hline
\end{tabular}
\end{center}
\end{table}

\begin{table}[ht!]
	\caption{\label{n20000}Results for $n=20,000$}
	\begin{center}
\begin{tabular}{|c||c|c|c|c|c|c|}						
	\hline						
	$p$&$\pi=0.0$&$\pi=0.2$&$\pi=0.4$&$\pi=0.6$&$\pi=0.8$&$\pi=1.0$\\						
	\hline						
	\multicolumn{7}{|c|}{Proportion of Market Share Captured (Power Decay)}\\						
	\hline						
	1&	0.11148&	0.11663&	0.12690&	0.13727&	0.14764&	0.15802\\
	2&	0.19983&	0.20479&	0.21919&	0.23453&	0.24986&	0.26520\\
	3&	0.27110&	0.27631&	0.28810&	0.30755&	0.32738&	0.34731\\
	4&	0.33213&	0.33692&	0.34882&	0.36732&	0.38808&	0.41197\\
	5&	0.38045&	0.38563&	0.39731&	0.41599&	0.43974&	0.46375\\
	10&	0.55228&	0.55426&	0.56319&	0.57819&	0.60184&	0.62784\\
	15&	0.65200&	0.65052&	0.65497&	0.66557&	0.68731&	0.71525\\
	20&	0.71637&	0.71531&	0.71717&	0.72304&	0.74114&	0.76902\\
	\hline\multicolumn{7}{|c|}{Run Time in Minutes for all 100 Runs (Power Decay)}\\\hline						
	1&	3.48&	6.05&	6.72&	14.20&	13.61&	14.35\\
	2&	7.07&	22.60&	25.50&	27.37&	26.66&	29.30\\
	3&	11.63&	30.98&	49.46&	46.79&	47.37&	51.29\\
	4&	14.70&	62.02&	69.76&	65.95&	66.74&	76.62\\
	5&	28.69&	64.37&	102.83&	95.19&	95.92&	104.75\\
	10&	87.67&	216.13&	262.15&	325.64&	326.74&	357.75\\
	15&	202.60&	418.13&	547.68&	599.07&	707.36&	816.69\\
	20&	377.48&	758.41&	849.22&	1006.90&	1280.02&	1524.14\\
	\hline						
	\multicolumn{7}{|c|}{Proportion of Market Share Captured (Exponential Decay)}\\						
	\hline						
	1&	0.12783&	0.12190&	0.14109&	0.16565&	0.19020&	0.21476\\
	2&	0.21591&	0.23018&	0.25634&	0.29040&	0.32447&	0.35854\\
	3&	0.28831&	0.30469&	0.33961&	0.38699&	0.43437&	0.48176\\
	4&	0.35422&	0.37689&	0.41490&	0.47210&	0.52942&	0.58675\\
	5&	0.40428&	0.42512&	0.46624&	0.52401&	0.58340&	0.64279\\
	10&	0.59344&	0.60850&	0.63644&	0.66946&	0.71626&	0.78028\\
	15&	0.69296&	0.70319&	0.71956&	0.74457&	0.78206&	0.84003\\
	20&	0.75125&	0.75820&	0.77252&	0.79319&	0.82232&	0.87325\\
	\hline						
	\multicolumn{7}{|c|}{Run Time in Minutes for all 100 Runs (Exponential Decay)}\\						
	\hline						
	1&	4.15&	11.90&	15.46&	15.10&	14.81&	14.93\\
	2&	8.57&	23.06&	31.68&	32.05&	31.62&	30.21\\
	3&	14.93&	42.87&	52.61&	67.63&	58.60&	56.91\\
	4&	23.80&	68.54&	78.36&	98.59&	91.84&	90.02\\
	5&	42.54&	94.38&	119.23&	141.04&	137.22&	140.87\\
	10&	137.63&	322.36&	421.90&	497.23&	519.75&	523.08\\
	15&	303.16&	743.76&	934.71&	1151.17&	1316.85&	1285.34\\
	20&	648.79&	1355.42&	1666.85&	1997.80&	2591.83&	2533.08\\
			\hline
		\end{tabular}
	\end{center}
\end{table}

\begin{table}[ht!]
	\caption{\label{Aver}Average Results for All Proportions for Each Instance}
	\begin{center}
				\setlength{\tabcolsep}{5pt}	
\begin{tabular}{|c||c|c|c|c|c|c|c|c|c|}									
	\hline									
	$p$&$n$=100&$n$=200&$n$=500&$n$=1000&$n$=2000&$n$=5000&$n$=10000&$n$=15000&$n$=20000\\									
	\hline									
	\multicolumn{10}{|c|}{Average Proportion of Market Share Captured (Power Decay)}\\									
	\hline									
	1&	0.12947&	0.13825&	0.13817&	0.13207&	0.13318&	0.13168&	0.13189&	0.13287&	0.13299\\
	2&	0.22831&	0.23300&	0.23492&	0.22843&	0.22950&	0.22703&	0.22712&	0.22849&	0.22890\\
	3&	0.30930&	0.30982&	0.31281&	0.30382&	0.30494&	0.30094&	0.30099&	0.30233&	0.30296\\
	4&	0.37307&	0.37121&	0.37441&	0.36445&	0.36600&	0.36163&	0.36161&	0.36331&	0.36421\\
	5&	0.42422&	0.42087&	0.42503&	0.41553&	0.41673&	0.41128&	0.41179&	0.41289&	0.41381\\
	10&	0.59155&	0.58731&	0.58877&	0.58119&	0.58279&	0.57736&	0.57786&	0.57897&	0.57960\\
	15&	0.68448&	0.68104&	0.68059&	0.67313&	0.67484&	0.67001&	0.66983&	0.67051&	0.67094\\
	20&	0.74347&	0.74012&	0.74069&	0.73343&	0.73460&	0.73050&	0.73001&	0.73035&	0.73034\\
	\hline\multicolumn{10}{|c|}{Average Run Time in Minutes for all 100 Runs (Power Decay)}\\\hline									
	1&	0.19&	0.29&	0.45&	0.73&	1.42&	2.92&	5.19&	7.39&	9.74\\
	2&	0.53&	0.52&	0.84&	1.49&	2.45&	5.78&	11.82&	17.48&	23.08\\
	3&	0.75&	0.91&	1.42&	2.21&	3.95&	9.64&	19.07&	28.37&	39.59\\
	4&	1.17&	1.22&	1.96&	3.21&	6.12&	15.38&	26.77&	41.61&	59.30\\
	5&	1.41&	1.61&	2.65&	4.44&	8.26&	21.53&	39.80&	58.63&	81.96\\
	10&	3.50&	4.47&	8.13&	13.28&	24.96&	64.97&	125.02&	188.70&	262.68\\
	15&	6.00&	8.40&	15.42&	27.49&	52.33&	132.67&	261.02&	406.75&	548.59\\
	20&	9.31&	13.83&	26.02&	47.21&	91.91&	236.82&	470.95&	753.85&	966.03\\
	\hline									
	\multicolumn{10}{|c|}{Average Proportion of Market Share Captured (Exponential Decay)}\\									
	\hline									
	1&	0.15753&	0.16695&	0.16806&	0.16151&	0.16236&	0.15875&	0.15841&	0.15978&	0.16024\\
	2&	0.28038&	0.28258&	0.29229&	0.28418&	0.28996&	0.28057&	0.27689&	0.27876&	0.27931\\
	3&	0.38560&	0.38153&	0.39516&	0.38054&	0.38123&	0.37367&	0.37150&	0.37202&	0.37262\\
	4&	0.46826&	0.45766&	0.46788&	0.45946&	0.46330&	0.45568&	0.45427&	0.45511&	0.45571\\
	5&	0.52203&	0.51386&	0.52170&	0.51129&	0.51429&	0.50597&	0.50548&	0.50668&	0.50764\\
	10&	0.69295&	0.67333&	0.67583&	0.66680&	0.66941&	0.66393&	0.66459&	0.66626&	0.66740\\
	15&	0.77460&	0.75354&	0.75409&	0.74605&	0.74909&	0.74558&	0.74557&	0.74636&	0.74706\\
	20&	0.82357&	0.80263&	0.80303&	0.79520&	0.79763&	0.79421&	0.79399&	0.79460&	0.79512\\
	\hline									
	\multicolumn{10}{|c|}{Average Run Time in Minutes for all 100 Runs (Exponential Decay)}\\									
	\hline									
	1&	0.34&	0.28&	0.42&	0.82&	1.38&	3.54&	6.36&	9.76&	12.73\\
	2&	0.66&	0.70&	1.13&	1.65&	2.94&	7.25&	13.48&	20.59&	26.20\\
	3&	1.04&	1.12&	1.86&	2.93&	5.10&	12.86&	24.11&	35.20&	48.93\\
	4&	1.45&	1.50&	2.61&	4.32&	8.25&	20.02&	37.75&	56.80&	75.19\\
	5&	1.69&	2.12&	4.14&	6.36&	12.15&	30.10&	55.23&	82.35&	112.55\\
	10&	4.84&	6.54&	13.13&	22.40&	43.65&	102.57&	203.63&	291.40&	403.66\\
	15&	10.26&	15.26&	29.94&	51.89&	96.61&	237.17&	471.89&	711.21&	955.83\\
	20&	18.75&	30.06&	55.13&	98.52&	185.07&	452.76&	895.43&	1324.70&	1798.96\\
			\hline
		\end{tabular}
	\end{center}
\end{table}

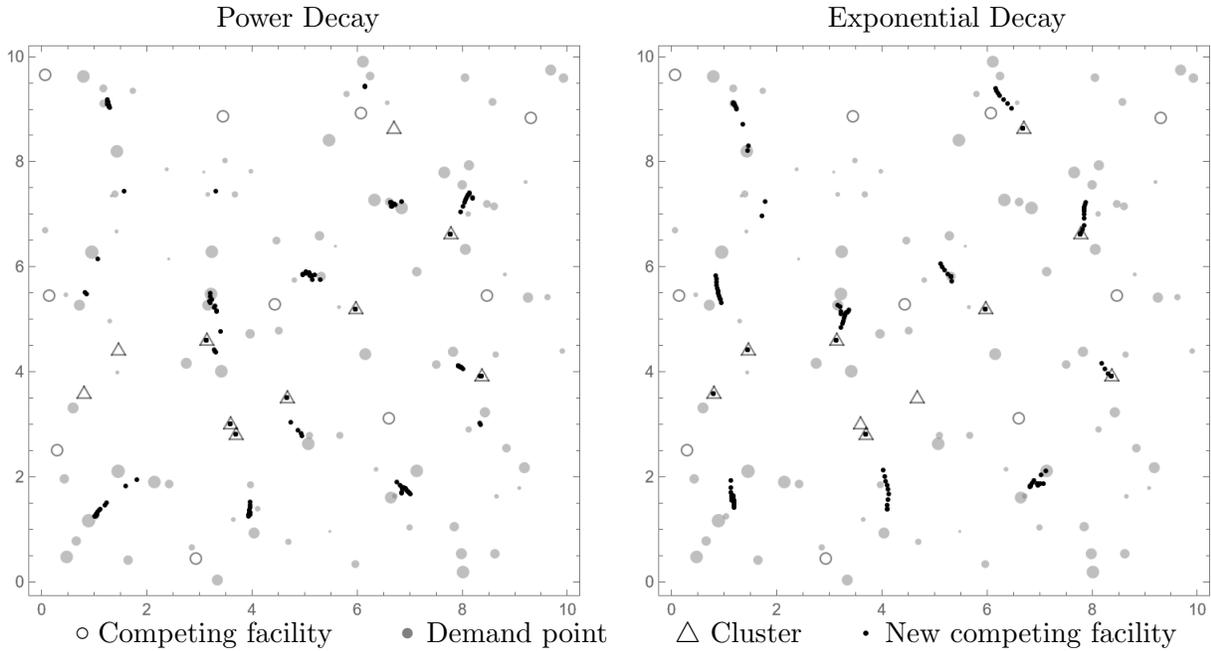
\begin{figure}[ht!]
	\centering
	\setlength{\unitlength}{1in}
	\begin{picture}(6.5,3.3)
	\includegraphics[width=3in]{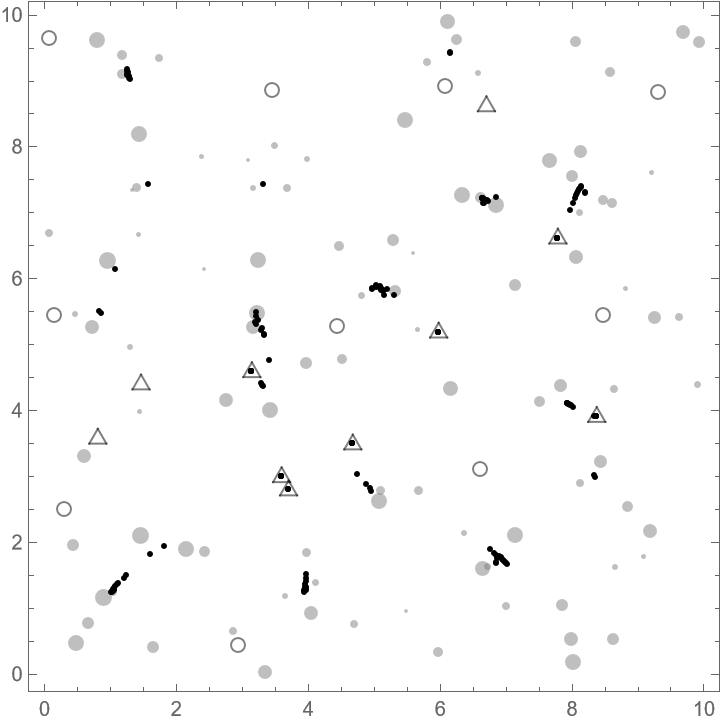}
	\hspace{0.2in}
	\includegraphics[width=3in]{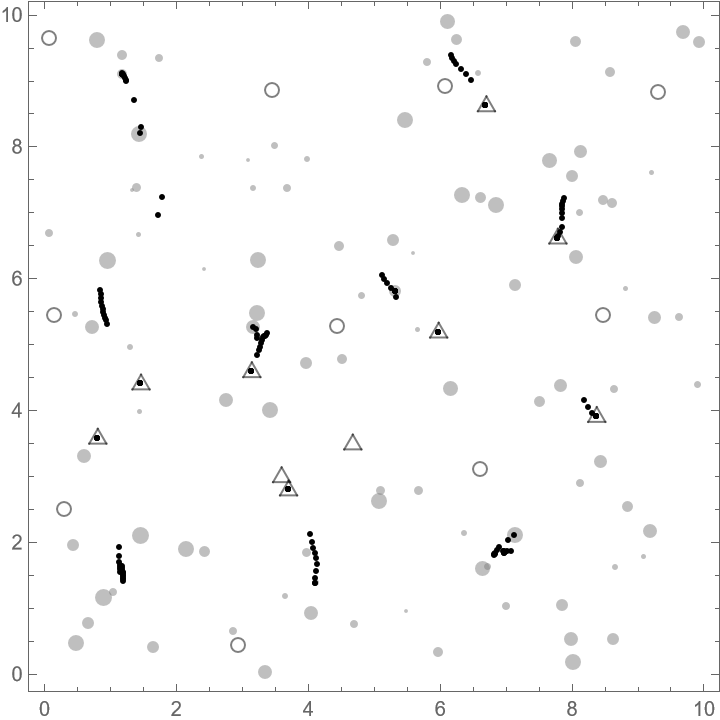}
	\put(-5.2,3.1){Power Decay}		
	\put(-2.,3.1){Exponential Decay}
\put(-2.8,-0.125){$\triangle$ Cluster}
\put(-4.2,-0.08){{{\color{gray}\circle*{0.06}}}}
\put(-4.1,-0.125){Demand point}
\put(-5.9,-0.08){{\circle{0.05}}}
\put(-5.82,-0.125){Competing facility}
\put(-1.8,-0.08){\circle*{0.03}}
\put(-1.7,-0.125){New competing facility}
	\end{picture}
	\caption{Distributed Locations of Facilities for $n=100$ and Various $\pi$ Values}
	\label{figloc}
\end{figure}

\begin{figure}[ht!]
	\centering
	\setlength{\unitlength}{1in}
	\begin{picture}(4,2.9)
	\includegraphics[width=4in]{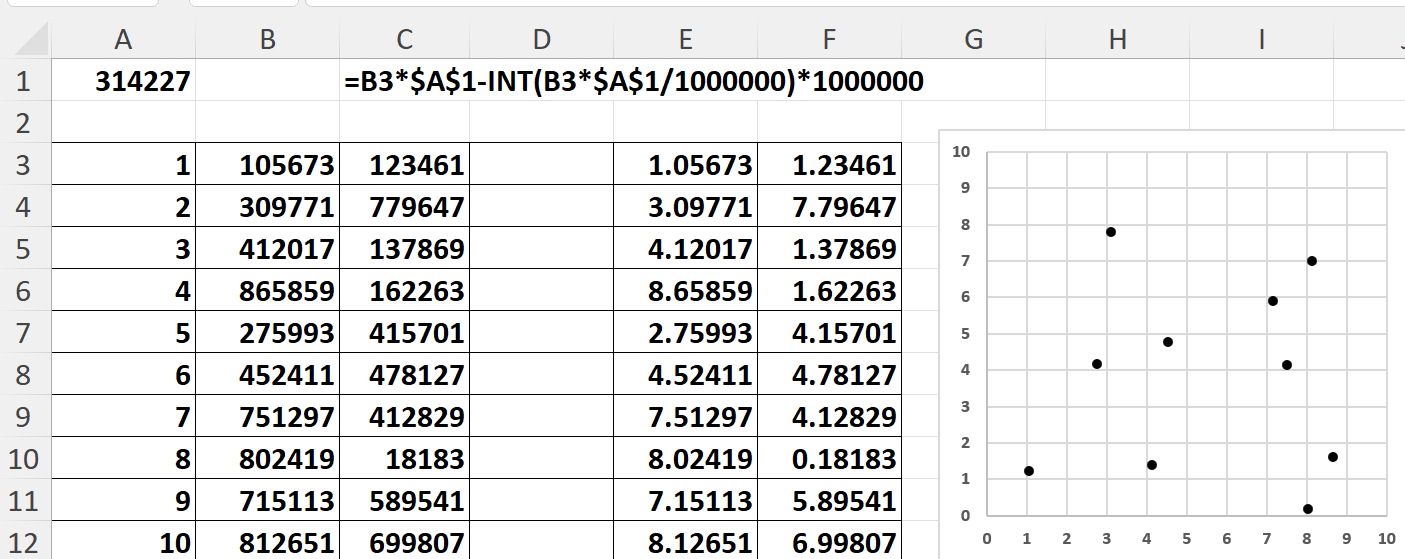}
    \put(-2.85,1.32){\vector(-1,-1){0.285in}}	
	\put(-2.08,1.23){\scriptsize $x$}
	\put(-1.65,1.24){\scriptsize $y$}
	\put(-4,2.3){\scriptsize The multiplier $\theta$ is in A1. In cells B3\&C3 put the $r_1$ seeds for the first point.}
	\put(-4,2.15){\scriptsize To get a different sequence, e.g. for facilities, just change the seeds in cells B3\&C3.}
	\put(-4,2){\scriptsize Copy B4 to C4. Then, B4\&C4 can be copied down up to 50,000 points.}
	\put(-4,1.85){\scriptsize $0\le x,y\le 10$ in columns E\&F are obtained by dividing columns B\&C by 100,000.}
	\put(-4,1.7){\scriptsize The first 10 points are plotted on the right.}
\thicklines
	\put(-4.1,-0.1){\line(1,0){4.5}}
	\put(0.4,-0.1){\line(0,1){2.6}}
	\put(-4.1,-0.1){\line(0,1){2.6}}
	\put(-4.1,2.5){\line(1,0){4.5}}
	\end{picture}
	\vspace{-0.2in}
	\caption{An Excel file to generate a pseudo-random sequence of numbers}
	\label{figex}
\end{figure}

\section{\label{sec3}Computational Experiments}

The experiments were conducted on a virtualized Linux Ubuntu 20.04 LTS environment with 32 vCPUs and 384GB of vRAM. The physical server used was a VxRail V570F with Intel Xeon Gold 6248 @ 2.5GHz (2 socket, 20 cores per processor, 80 logical processors), 748GB RAM, and VMware vSAN storage.
 For power decay we applied $\lambda=2$, and for exponential decay $\lambda=1$.

In the computational experiments we tested only problems for which all $p$ facilities are new and there are no existing chain facilities. We considered the problems investigated in \citet{ZOK23,DOD22}. There are 10 existing competing facilities and 10 clusters. These instances, that can be easily replicated for future comparisons, are generated by a pseudo random number generator based on the one proposed by \citet{LK91}.
Integer starting seed $1\le r_1\le 999,999$ and an odd multiplier $\theta$, which are not divisible by 5, are selected. We used $\theta=314,227$. The sequence is generated by the following rule for $k\ge 1$:
$
r_{k+1}=\theta r_k-\lfloor\frac{\theta r_k}{1,000,000}\rfloor\times 1,000,000.
$
To generate a sequence of numbers in the range $a<x_k<b$, $x_k =a+(b-a)\frac{r_k}{1,000,000}$ is calculated.  A simple Excel file that can produce such sequences is depicted in Figure \ref{figex}. On the right side of the figure the first 10 points are depicted.
 For complete details, and a FORTRAN code for generating these instances, see \citet{DOD22}.
 
The instance  with $n=100$ demand points is depicted in Figure \ref{figloc}.
The results for $n=100$ and 20000 are reported in Tables \ref{n100}-\ref{n20000}. The results for $n=200, 500, 1000, 2000, 5000, 10000$, and 15000 are reported in the appendix. The average performance for all values of $p$ and $n$ instances are summarized in Table \ref{Aver}.

	As reported in Table \ref{Aver}, average run times (for all 100 starts) were about 13 minutes for the largest instance of $p=1$, and increased to an average of 1800 minutes for $p=20$. By a regression analysis on the average run times we found for $p=20$ that: (i) for power decay run times are approximately proportional to $n^{0.905}$ with a p-value of $1.3\cdot 10^{-8}$; and (ii) for exponential decay they are proportional to $n^{0.875}$ with a p-value of $3.8\cdot 10^{-9}$. This means that run times are less than linear in the number of demand points. A regression analysis on all 72 values for each distance decay function in Table \ref{Aver} found that run time is proportional to $n^{0.815}p^{1.436}$ with a p-value of $2.1\cdot10^{-62}$ for power decay, and $n^{0.812}p^{1.613}$ with a p-value of $6.5\cdot10^{-61}$ for exponential decay. Such performance shows that the solution approach is very effective.


\subsection{Locations as a Function of $\pi$}

\citet{DOD22}, who investigated the location of one new competing facility,  depicted (in Figure 2 of that paper) the trajectory of the locations of the new  competing facility for $0\le \pi\le 1$ and $n=100$ demand points. The trajectory is discontinuous. Up to $\pi=0.28$ the locations are close to one another, at $\pi=0.29$ they ``jump" to another region and stay there until $\pi=0.79$, and so on. The locations of $p=10$ facilities among 100 demand points, investigated in the current paper, are depicted in Figure~\ref{figloc}. There is no clear trajectory for each facility. They are discontinuous as well. For multiple facilities there is no clear determination which ``jump" belongs to which facility. There is no meaning to the order of the facilities. They all have the same attractiveness so they are interchangeable leading to $p!$ trajectories.

One phenomenon, that confirms the observation in previous papers, is that many facilities tend to locate at clusters. For power decay, there are seven clusters that have facilities located at them for at least one value of $\pi$. These are the triangles with a black dot inside. For exponential decay there are eight such clusters. For example, consider the case that the competing facilities are coffee shops like Starbucks, and the clusters are grocery stores. We observed in our neighborhood, which is typical, that a significant proportion of coffee shops tend to be located in grocery stores. These are good locations because many customers get coffee in conjunction with a trip to a grocery store.
		
\subsection{Discussion of the Results}
All results for the instances reported in Tables \ref{n100}-\ref{n20000} and Tables \ref{n200}-\ref{n15000} in the Appendix, have a common structure. The market share increases as $\pi$ increases and the largest market share captured by the new facilities is the highest for $\pi=1$, which is when all trips are multipurpose trips. As can be expected, the expected market share captured increases with an increase in the number of new competing facilities as reported in Table \ref{Aver}. Adding more facilities cannibalizes some of the market share captured by other facilities but a better configuration can be established and the additional market share attracted from existing competing facilities exceeds the market share lost by cannibalization \citep{MMR20,Dr10,Pl05}. In the tested instances, there are 10 existing facilities and $p$ new ones. If the new facilities are randomly located, we would expect a proportion of $\frac{p}{10+p}$ of the total buying power to be attracted to the new facilities. This entails $\frac{1}{11}$=0.091 for $p=1$, $\frac{10}{20}$=0.500 for $p=10$, and $\frac{20}{30}$=0.667 for $p=20$. As should be expected, the averages reported in Table \ref{Aver} exceed these expected captured market shares because locations which are better than random are found by the solution procedures. We calculated the ratio between the expected market share for the average of all $n$'s in Table \ref{Aver}, with the expected market share for random locations. The ratios are 1.47 for $p=1$ decreasing as $p$ increases to a ratio of 1.10 for $p=20$ (for larger values of $p$ there is less room for an increase). For exponential decay the ratios are 1.78 decreasing to 1.20. Exponential decay is more sensitive to the locations of the new facilities.

The calculated market share captured by the exponential decay function is greater than the market share captured by the power decay function (Table \ref{Aver}). This is a good result because  \citet{Dr06a} showed that exponential decay estimates the captured market share better than power decay for shopping malls.

Run times generally increase with the increase of $\pi$. It usually requires more time to find a solution when a larger proportion of the trips are multipurpose. Exponential decay run times are higher than power decay.

\section{\label{sec4}Conclusions and Suggestions for Future Research}

This paper continues to develop and explore the impact of multipurpose (MP) trips on retail location. While there has been longstanding concern for this topic, we develop new insights by allowing the site selector to pick several sites. The addition of  multiple new outlets can cause cannibalization of existing sales, but this effect can be mitigated by careful arrangement of the selections. We also allow exploration of alternative functional forms for the distance decay effect (power and exponential). In keeping with past results, the exponential formulation delivers greater capture of market share. Generally, the introduction of MP trips provides a better estimate of the total market share and consequently the location decision is more accurate as well. 

\subsection{Suggestions for Future Research}

The same models can be applied for different estimation of the captured market share rather than the gravity model.
Other useful options to apply and test include (we expect shorter run times for the same $n$ and $p$):

	\begin{enumerate}
			\item We can modify  the model with the additional option that some of the competitors are moved or closed.
		\item To consider closing or replacing some of our stores in an existing market, possibly with the additional possibility that we add new replacements. To find the optimal solution it requires testing the final market share when all possible $q$ facilities are removed and $p$ facilities added which is applying our approach for every selection of $q$ facilities. This is practical for $q=1$. Since closing facilities can be very costly, it may suffice to test the possibility of closing only one facility, and hope that the increase in market share will exceed the closing cost.
	\item	\citet{DZ23} considered the possibility that one of the competing facilities goes out of business, and their approach can be modified and implemented for this case.
\item Locating one new facility when there are several existing facilities that belong to the same chain can be optimally solved by the BSSS \citep{HPT81} or BTST \citep{DS04} global optimization algorithms, which find the optimal solution within a given relative accuracy.
	\end{enumerate}

\renewcommand{\baselinestretch}{1}
\renewcommand{\arraystretch}{1}
\large
\normalsize

	
	\bibliographystyle{apalike}

\renewcommand{\baselinestretch}{2}
\renewcommand{\arraystretch}{0.5}
\large
\normalsize

\newpage\clearpage

{\centering{\section*{Appendix}}}
\begin{table}[ht!]
	\caption{\label{n200}Results for $n=200$}
	\begin{center}
\begin{tabular}{|c||c|c|c|c|c|c|}						
	\hline						
	$p$&$\pi=0.0$&$\pi=0.2$&$\pi=0.4$&$\pi=0.6$&$\pi=0.8$&$\pi=1.0$\\						
	\hline						
	\multicolumn{7}{|c|}{Proportion of Market Share Captured (Power Decay)}\\						
	\hline						
	1&	0.12275&	0.12427&	0.13246&	0.14124&	0.15001&	0.15879\\
	2&	0.21513&	0.21637&	0.22304&	0.23310&	0.24679&	0.26359\\
	3&	0.29254&	0.29387&	0.29595&	0.30899&	0.32525&	0.34232\\
	4&	0.35228&	0.35292&	0.35794&	0.36784&	0.38794&	0.40834\\
	5&	0.40467&	0.40185&	0.40456&	0.41534&	0.43759&	0.46119\\
	10&	0.58303&	0.57474&	0.56875&	0.57566&	0.59714&	0.62453\\
	15&	0.68543&	0.67386&	0.66610&	0.66592&	0.68352&	0.71144\\
	20&	0.74601&	0.73564&	0.72847&	0.72587&	0.73894&	0.76577\\
	\hline\multicolumn{7}{|c|}{Run Time in Minutes for all 100 Runs (Power Decay)}\\\hline						
	1&	0.10&	0.09&	0.33&	0.40&	0.40&	0.40\\
	2&	0.22&	0.24&	0.36&	0.78&	0.75&	0.77\\
	3&	0.34&	0.41&	0.80&	1.20&	1.26&	1.48\\
	4&	0.48&	0.55&	1.20&	1.66&	1.71&	1.76\\
	5&	0.63&	0.78&	1.38&	2.20&	2.27&	2.39\\
	10&	1.79&	1.95&	3.22&	5.25&	7.44&	7.20\\
	15&	3.30&	3.28&	5.38&	9.57&	13.71&	15.17\\
	20&	5.42&	5.55&	8.32&	15.62&	21.75&	26.30\\
	\hline						
	\multicolumn{7}{|c|}{Proportion of Market Share Captured (Exponential Decay)}\\						
	\hline						
	1&	0.13804&	0.13010&	0.14339&	0.17006&	0.19672&	0.22339\\
	2&	0.23168&	0.23053&	0.25080&	0.28895&	0.32711&	0.36642\\
	3&	0.30950&	0.32890&	0.35210&	0.38353&	0.43264&	0.48254\\
	4&	0.37585&	0.39203&	0.41932&	0.45793&	0.51941&	0.58144\\
	5&	0.43379&	0.44206&	0.47290&	0.51866&	0.57564&	0.64008\\
	10&	0.62167&	0.62510&	0.64298&	0.66681&	0.70898&	0.77444\\
	15&	0.71734&	0.71607&	0.72761&	0.74800&	0.77769&	0.83450\\
	20&	0.77579&	0.77305&	0.78074&	0.79541&	0.82029&	0.87050\\
	\hline						
	\multicolumn{7}{|c|}{Run Time in Minutes for all 100 Runs (Exponential Decay)}\\						
	\hline						
	1&	0.10&	0.11&	0.19&	0.46&	0.49&	0.35\\
	2&	0.23&	0.59&	0.77&	0.90&	0.83&	0.88\\
	3&	0.45&	0.89&	1.08&	1.51&	1.40&	1.39\\
	4&	0.74&	0.75&	1.44&	2.01&	2.04&	2.00\\
	5&	1.21&	1.11&	2.27&	2.52&	2.90&	2.74\\
	10&	3.64&	4.69&	5.46&	7.96&	8.56&	8.92\\
	15&	10.57&	10.71&	14.54&	17.11&	18.74&	19.87\\
	20&	22.52&	29.62&	28.84&	29.57&	32.63&	37.15\\
			\hline
		\end{tabular}
	\end{center}
\end{table}

\begin{table}[ht!]
	\caption{\label{n500}Results for $n=500$}
	\begin{center}
\begin{tabular}{|c||c|c|c|c|c|c|}						
	\hline						
	$p$&$\pi=0.0$&$\pi=0.2$&$\pi=0.4$&$\pi=0.6$&$\pi=0.8$&$\pi=1.0$\\						
	\hline						
	\multicolumn{7}{|c|}{Proportion of Market Share Captured (Power Decay)}\\						
	\hline						
	1&	0.12185&	0.12357&	0.13236&	0.14139&	0.15041&	0.15944\\
	2&	0.21279&	0.21481&	0.22464&	0.23854&	0.25243&	0.26632\\
	3&	0.29013&	0.29243&	0.29910&	0.31362&	0.33146&	0.35010\\
	4&	0.34968&	0.35341&	0.36077&	0.37528&	0.39316&	0.41416\\
	5&	0.39947&	0.40280&	0.41050&	0.42473&	0.44572&	0.46698\\
	10&	0.56810&	0.56764&	0.57253&	0.58623&	0.60718&	0.63092\\
	15&	0.67024&	0.66578&	0.66472&	0.67195&	0.69307&	0.71779\\
	20&	0.73777&	0.73129&	0.72736&	0.72978&	0.74659&	0.77136\\
	\hline\multicolumn{7}{|c|}{Run Time in Minutes for all 100 Runs (Power Decay)}\\\hline						
	1&	0.20&	0.21&	0.57&	0.54&	0.58&	0.60\\
	2&	0.31&	0.31&	0.98&	1.13&	1.09&	1.21\\
	3&	0.52&	0.73&	1.76&	1.75&	1.83&	1.91\\
	4&	0.77&	0.94&	2.19&	2.46&	2.59&	2.81\\
	5&	1.08&	1.33&	2.36&	3.48&	3.67&	3.97\\
	10&	3.06&	4.17&	6.91&	11.01&	11.32&	12.32\\
	15&	6.04&	7.59&	11.83&	18.44&	23.00&	25.60\\
	20&	10.56&	12.77&	19.04&	28.84&	39.91&	45.00\\
	\hline						
	\multicolumn{7}{|c|}{Proportion of Market Share Captured (Exponential Decay)}\\						
	\hline						
	1&	0.14272&	0.13620&	0.14308&	0.16921&	0.19545&	0.22173\\
	2&	0.23140&	0.24197&	0.26539&	0.30185&	0.33832&	0.37478\\
	3&	0.31535&	0.33253&	0.36160&	0.40517&	0.45383&	0.50250\\
	4&	0.37707&	0.39601&	0.42741&	0.47870&	0.53559&	0.59248\\
	5&	0.42500&	0.44292&	0.48170&	0.53483&	0.59352&	0.65221\\
	10&	0.60308&	0.61746&	0.64323&	0.67774&	0.72641&	0.78705\\
	15&	0.70290&	0.70999&	0.72592&	0.75028&	0.78948&	0.84595\\
	20&	0.76633&	0.76846&	0.77925&	0.79846&	0.82755&	0.87812\\
	\hline						
	\multicolumn{7}{|c|}{Run Time in Minutes for all 100 Runs (Exponential Decay)}\\						
	\hline						
	1&	0.15&	0.32&	0.32&	0.42&	0.67&	0.66\\
	2&	0.32&	1.07&	1.68&	1.25&	1.26&	1.23\\
	3&	0.57&	1.64&	2.31&	2.19&	2.22&	2.21\\
	4&	0.94&	1.37&	3.03&	3.56&	3.34&	3.39\\
	5&	1.58&	2.83&	5.65&	5.27&	4.83&	4.67\\
	10&	5.84&	10.18&	12.96&	16.01&	17.15&	16.67\\
	15&	13.41&	26.10&	27.06&	35.17&	39.83&	38.08\\
	20&	32.65&	40.27&	48.22&	62.90&	74.32&	72.43\\
			\hline
		\end{tabular}
	\end{center}
\end{table}

\begin{table}[ht!]
	\caption{\label{n1000}Results for $n=1,000$}
	\begin{center}
\begin{tabular}{|c||c|c|c|c|c|c|}						
	\hline						
	$p$&$\pi=0.0$&$\pi=0.2$&$\pi=0.4$&$\pi=0.6$&$\pi=0.8$&$\pi=1.0$\\						
	\hline						
	\multicolumn{7}{|c|}{Proportion of Market Share Captured (Power Decay)}\\						
	\hline						
	1&	0.11853&	0.11578&	0.12476&	0.13460&	0.14444&	0.15429\\
	2&	0.20830&	0.21004&	0.21512&	0.23015&	0.24570&	0.26126\\
	3&	0.27969&	0.28183&	0.28846&	0.30493&	0.32422&	0.34380\\
	4&	0.33791&	0.34084&	0.34920&	0.36485&	0.38512&	0.40877\\
	5&	0.38944&	0.39073&	0.39930&	0.41468&	0.43781&	0.46124\\
	10&	0.56025&	0.55869&	0.56333&	0.57831&	0.60056&	0.62597\\
	15&	0.65992&	0.65613&	0.65666&	0.66561&	0.68670&	0.71377\\
	20&	0.72618&	0.72160&	0.71995&	0.72389&	0.74100&	0.76793\\
	\hline\multicolumn{7}{|c|}{Run Time in Minutes for all 100 Runs (Power Decay)}\\\hline						
	1&	0.21&	0.47&	0.80&	0.75&	1.06&	1.08\\
	2&	0.47&	1.29&	1.74&	1.77&	1.83&	1.86\\
	3&	0.79&	0.88&	2.71&	2.86&	2.94&	3.05\\
	4&	1.20&	1.77&	3.31&	4.27&	4.33&	4.39\\
	5&	1.58&	2.51&	4.43&	5.78&	6.20&	6.10\\
	10&	4.92&	5.34&	11.88&	16.90&	20.29&	20.38\\
	15&	11.05&	13.43&	23.88&	31.42&	41.78&	43.41\\
	20&	20.00&	25.14&	38.04&	49.46&	71.37&	79.27\\
	\hline						
	\multicolumn{7}{|c|}{Proportion of Market Share Captured (Exponential Decay)}\\						
	\hline						
	1&	0.13897&	0.13224&	0.13735&	0.16209&	0.18683&	0.21157\\
	2&	0.22955&	0.23879&	0.25645&	0.29161&	0.32677&	0.36193\\
	3&	0.29771&	0.31559&	0.34472&	0.39278&	0.44174&	0.49069\\
	4&	0.36463&	0.38399&	0.41772&	0.47270&	0.53014&	0.58758\\
	5&	0.41549&	0.43304&	0.47107&	0.52298&	0.58272&	0.64246\\
	10&	0.59738&	0.60539&	0.63088&	0.66927&	0.71807&	0.77981\\
	15&	0.69408&	0.70116&	0.71763&	0.74170&	0.78179&	0.83993\\
	20&	0.75422&	0.75913&	0.77231&	0.79162&	0.82031&	0.87358\\
	\hline						
	\multicolumn{7}{|c|}{Run Time in Minutes for all 100 Runs (Exponential Decay)}\\						
	\hline						
	1&	0.20&	0.69&	0.90&	1.09&	1.01&	1.04\\
	2&	0.53&	1.60&	1.75&	2.09&	2.00&	1.94\\
	3&	0.94&	2.09&	3.30&	3.87&	3.87&	3.49\\
	4&	1.47&	3.09&	4.57&	6.00&	5.53&	5.23\\
	5&	2.48&	4.33&	6.70&	8.51&	8.49&	7.66\\
	10&	7.88&	18.28&	21.99&	26.35&	31.45&	28.43\\
	15&	21.80&	35.56&	46.99&	62.56&	76.57&	67.88\\
	20&	41.16&	74.95&	92.11&	113.48&	138.44&	130.96\\
			\hline
		\end{tabular}
	\end{center}
\end{table}

\begin{table}[ht!]
	\caption{\label{n2000}Results for $n=2,000$}
	\begin{center}
\begin{tabular}{|c||c|c|c|c|c|c|}						
	\hline						
	$p$&$\pi=0.0$&$\pi=0.2$&$\pi=0.4$&$\pi=0.6$&$\pi=0.8$&$\pi=1.0$\\						
	\hline						
	\multicolumn{7}{|c|}{Proportion of Market Share Captured (Power Decay)}\\						
	\hline						
	1&	0.11613&	0.11735&	0.12693&	0.13652&	0.14610&	0.15603\\
	2&	0.20742&	0.21053&	0.21750&	0.23164&	0.24713&	0.26276\\
	3&	0.27746&	0.28117&	0.29267&	0.30749&	0.32549&	0.34535\\
	4&	0.33667&	0.34133&	0.35179&	0.36872&	0.38697&	0.41054\\
	5&	0.38756&	0.39131&	0.40047&	0.41889&	0.43939&	0.46275\\
	10&	0.56085&	0.55951&	0.56562&	0.58054&	0.60293&	0.62728\\
	15&	0.66194&	0.65765&	0.65817&	0.66783&	0.68844&	0.71498\\
	20&	0.72755&	0.72207&	0.72112&	0.72542&	0.74215&	0.76932\\
	\hline\multicolumn{7}{|c|}{Run Time in Minutes for all 100 Runs (Power Decay)}\\\hline						
	1&	0.42&	1.64&	1.41&	1.27&	2.13&	1.64\\
	2&	0.84&	0.95&	3.19&	3.21&	3.11&	3.38\\
	3&	1.44&	1.53&	5.08&	5.14&	5.04&	5.45\\
	4&	1.99&	3.82&	7.52&	7.99&	7.65&	7.78\\
	5&	2.91&	5.48&	8.52&	10.93&	10.41&	11.27\\
	10&	9.31&	12.23&	25.33&	32.59&	34.52&	35.78\\
	15&	19.77&	28.21&	44.01&	62.93&	77.09&	81.97\\
	20&	37.09&	57.89&	78.63&	96.21&	128.83&	152.82\\
	\hline						
	\multicolumn{7}{|c|}{Proportion of Market Share Captured (Exponential Decay)}\\						
	\hline						
	1&	0.13323&	0.12828&	0.14100&	0.16577&	0.19054&	0.21532\\
	2&	0.22616&	0.23956&	0.26586&	0.30096&	0.33606&	0.37115\\
	3&	0.29389&	0.31510&	0.34760&	0.39560&	0.44360&	0.49160\\
	4&	0.35957&	0.38550&	0.42358&	0.47996&	0.53705&	0.59415\\
	5&	0.41097&	0.43229&	0.47388&	0.52989&	0.58955&	0.64920\\
	10&	0.59766&	0.60882&	0.63641&	0.67090&	0.71964&	0.78302\\
	15&	0.69598&	0.70535&	0.72128&	0.74543&	0.78409&	0.84238\\
	20&	0.75674&	0.76142&	0.77409&	0.79438&	0.82331&	0.87585\\
	\hline						
	\multicolumn{7}{|c|}{Run Time in Minutes for all 100 Runs (Exponential Decay)}\\						
	\hline						
	1&	0.37&	1.11&	1.68&	1.67&	1.82&	1.65\\
	2&	0.92&	2.54&	3.96&	3.54&	3.43&	3.23\\
	3&	1.71&	3.21&	6.16&	6.78&	6.43&	6.30\\
	4&	2.95&	7.45&	8.70&	10.36&	10.23&	9.80\\
	5&	4.33&	8.82&	14.58&	15.26&	16.02&	13.89\\
	10&	18.53&	39.38&	45.34&	49.65&	55.62&	53.37\\
	15&	36.59&	69.54&	89.60&	117.14&	136.44&	130.35\\
	20&	66.72&	139.56&	171.96&	207.15&	264.26&	260.74\\
			\hline
		\end{tabular}
	\end{center}
\end{table}

\begin{table}[ht!]
	\caption{\label{n5000}Results for $n=5,000$}
	\begin{center}
\begin{tabular}{|c||c|c|c|c|c|c|}						
	\hline						
	$p$&$\pi=0.0$&$\pi=0.2$&$\pi=0.4$&$\pi=0.6$&$\pi=0.8$&$\pi=1.0$\\						
	\hline						
	\multicolumn{7}{|c|}{Proportion of Market Share Captured (Power Decay)}\\						
	\hline						
	1&	0.11125&	0.11566&	0.12535&	0.13556&	0.14594&	0.15631\\
	2&	0.20047&	0.20565&	0.21568&	0.23117&	0.24680&	0.26242\\
	3&	0.26862&	0.27473&	0.28879&	0.30434&	0.32447&	0.34471\\
	4&	0.32686&	0.33360&	0.34797&	0.36650&	0.38542&	0.40941\\
	5&	0.37640&	0.38176&	0.39651&	0.41501&	0.43685&	0.46113\\
	10&	0.55086&	0.55166&	0.55945&	0.57606&	0.60044&	0.62572\\
	15&	0.65289&	0.65061&	0.65352&	0.66338&	0.68573&	0.71395\\
	20&	0.71951&	0.71618&	0.71717&	0.72252&	0.73999&	0.76760\\
	\hline\multicolumn{7}{|c|}{Run Time in Minutes for all 100 Runs (Power Decay)}\\\hline						
	1&	1.00&	3.59&	2.17&	3.29&	3.79&	3.69\\
	2&	1.77&	3.30&	7.60&	7.29&	7.42&	7.29\\
	3&	3.00&	5.83&	12.34&	11.87&	12.29&	12.49\\
	4&	5.93&	10.65&	20.19&	19.03&	17.56&	18.92\\
	5&	7.80&	15.98&	29.17&	25.87&	24.91&	25.46\\
	10&	22.80&	45.81&	68.98&	82.10&	82.19&	87.93\\
	15&	50.88&	90.92&	128.74&	155.24&	175.93&	194.33\\
	20&	96.05&	161.09&	214.41&	251.30&	340.52&	357.56\\
	\hline						
	\multicolumn{7}{|c|}{Proportion of Market Share Captured (Exponential Decay)}\\						
	\hline						
	1&	0.12726&	0.12228&	0.13930&	0.16359&	0.18788&	0.21217\\
	2&	0.21493&	0.22962&	0.25886&	0.29277&	0.32667&	0.36057\\
	3&	0.28222&	0.30207&	0.34300&	0.39062&	0.43824&	0.48586\\
	4&	0.34565&	0.37223&	0.41773&	0.47526&	0.53283&	0.59040\\
	5&	0.39737&	0.42066&	0.46576&	0.52412&	0.58384&	0.64406\\
	10&	0.58816&	0.60106&	0.63005&	0.66748&	0.71627&	0.78055\\
	15&	0.69124&	0.70097&	0.71762&	0.74160&	0.78193&	0.84011\\
	20&	0.75048&	0.75666&	0.77137&	0.79225&	0.82061&	0.87391\\
	\hline						
	\multicolumn{7}{|c|}{Run Time in Minutes for all 100 Runs (Exponential Decay)}\\						
	\hline						
	1&	1.18&	3.59&	4.28&	4.19&	3.72&	4.30\\
	2&	2.43&	6.03&	9.45&	9.79&	7.70&	8.09\\
	3&	4.26&	11.17&	16.04&	15.96&	14.69&	15.08\\
	4&	7.69&	17.13&	24.12&	25.74&	22.57&	22.89\\
	5&	10.38&	23.83&	32.28&	44.14&	35.00&	34.99\\
	10&	34.36&	76.98&	107.59&	123.24&	137.91&	135.36\\
	15&	78.84&	171.62&	231.25&	280.55&	342.31&	318.49\\
	20&	158.73&	345.90&	429.98&	504.34&	640.17&	637.46\\
			\hline
		\end{tabular}
	\end{center}
\end{table}

\begin{table}[ht!]
	\caption{\label{n10000}Results for $n=10,000$}
	\begin{center}
\begin{tabular}{|c||c|c|c|c|c|c|}						
	\hline						
	$p$&$\pi=0.0$&$\pi=0.2$&$\pi=0.4$&$\pi=0.6$&$\pi=0.8$&$\pi=1.0$\\						
	\hline						
	\multicolumn{7}{|c|}{Proportion of Market Share Captured (Power Decay)}\\						
	\hline						
	1&	0.10906&	0.11588&	0.12612&	0.13644&	0.14677&	0.15709\\
	2&	0.19752&	0.20333&	0.21729&	0.23274&	0.24819&	0.26364\\
	3&	0.26858&	0.27399&	0.28665&	0.30554&	0.32553&	0.34568\\
	4&	0.32732&	0.33325&	0.34700&	0.36570&	0.38612&	0.41024\\
	5&	0.37735&	0.38216&	0.39606&	0.41473&	0.43823&	0.46219\\
	10&	0.54973&	0.55099&	0.56134&	0.57703&	0.60112&	0.62697\\
	15&	0.65133&	0.64937&	0.65338&	0.66431&	0.68644&	0.71417\\
	20&	0.71725&	0.71482&	0.71662&	0.72240&	0.74042&	0.76854\\
	\hline\multicolumn{7}{|c|}{Run Time in Minutes for all 100 Runs (Power Decay)}\\\hline						
	1&	1.49&	4.00&	3.91&	7.34&	7.12&	7.30\\
	2&	3.27&	12.88&	14.22&	13.26&	13.37&	13.92\\
	3&	5.69&	16.43&	24.38&	22.45&	23.29&	22.17\\
	4&	7.66&	17.69&	31.96&	35.76&	33.51&	34.08\\
	5&	12.78&	30.92&	53.15&	49.18&	46.72&	46.07\\
	10&	43.05&	95.89&	134.21&	158.65&	158.34&	160.01\\
	15&	108.60&	184.86&	246.22&	301.15&	357.36&	367.94\\
	20&	179.33&	365.30&	428.31&	504.33&	638.40&	710.02\\
	\hline						
	\multicolumn{7}{|c|}{Proportion of Market Share Captured (Exponential Decay)}\\						
	\hline						
	1&	0.12390&	0.11863&	0.14073&	0.16490&	0.18907&	0.21324\\
	2&	0.21069&	0.22667&	0.25564&	0.28921&	0.32277&	0.35634\\
	3&	0.28089&	0.30135&	0.34085&	0.38807&	0.43530&	0.48252\\
	4&	0.34736&	0.37189&	0.41570&	0.47296&	0.53023&	0.58750\\
	5&	0.39961&	0.42025&	0.46521&	0.52323&	0.58249&	0.64209\\
	10&	0.58800&	0.60366&	0.63196&	0.66756&	0.71598&	0.78039\\
	15&	0.69023&	0.70058&	0.71764&	0.74270&	0.78197&	0.84030\\
	20&	0.74962&	0.75634&	0.77104&	0.79202&	0.82165&	0.87330\\
	\hline						
	\multicolumn{7}{|c|}{Run Time in Minutes for all 100 Runs (Exponential Decay)}\\						
	\hline						
	1&	2.05&	6.41&	7.97&	7.89&	7.16&	6.71\\
	2&	4.77&	11.57&	15.68&	17.72&	15.44&	15.69\\
	3&	7.44&	20.57&	28.89&	30.99&	29.73&	27.05\\
	4&	12.40&	33.89&	46.06&	48.11&	44.11&	41.95\\
	5&	22.70&	43.12&	61.18&	74.41&	66.45&	63.50\\
	10&	80.11&	161.94&	216.16&	245.22&	262.33&	256.03\\
	15&	174.75&	357.00&	468.95&	564.75&	649.61&	616.29\\
	20&	322.03&	721.49&	815.82&	996.31&	1248.11&	1268.80\\
			\hline
		\end{tabular}
	\end{center}
\end{table}

\begin{table}[ht!]
	\caption{\label{n15000}Results for $n=15,000$}
	\begin{center}
\begin{tabular}{|c||c|c|c|c|c|c|}						
	\hline						
	$p$&$\pi=0.0$&$\pi=0.2$&$\pi=0.4$&$\pi=0.6$&$\pi=0.8$&$\pi=1.0$\\						
	\hline						
	\multicolumn{7}{|c|}{Proportion of Market Share Captured (Power Decay)}\\						
	\hline						
	1&	0.11036&	0.11676&	0.12701&	0.13736&	0.14771&	0.15806\\
	2&	0.19891&	0.20443&	0.21883&	0.23421&	0.24959&	0.26497\\
	3&	0.26950&	0.27512&	0.28793&	0.30719&	0.32711&	0.34713\\
	4&	0.32947&	0.33487&	0.34868&	0.36724&	0.38780&	0.41178\\
	5&	0.37843&	0.38381&	0.39699&	0.41548&	0.43915&	0.46348\\
	10&	0.55125&	0.55335&	0.56233&	0.57786&	0.60170&	0.62734\\
	15&	0.65184&	0.64985&	0.65444&	0.66513&	0.68704&	0.71475\\
	20&	0.71669&	0.71542&	0.71727&	0.72284&	0.74094&	0.76894\\
	\hline\multicolumn{7}{|c|}{Run Time in Minutes for all 100 Runs (Power Decay)}\\\hline						
	1&	2.74&	4.30&	6.30&	10.75&	9.64&	10.60\\
	2&	5.08&	18.52&	22.30&	18.86&	18.67&	21.43\\
	3&	8.45&	23.70&	38.47&	32.72&	31.71&	35.19\\
	4&	11.37&	32.48&	55.86&	50.22&	46.58&	53.16\\
	5&	23.27&	41.83&	76.11&	69.53&	68.44&	72.58\\
	10&	63.06&	165.26&	207.53&	227.89&	227.90&	240.54\\
	15&	173.68&	307.40&	403.64&	440.91&	514.73&	600.12\\
	20&	287.81&	719.51&	681.77&	740.03&	974.28&	1119.69\\
	\hline						
	\multicolumn{7}{|c|}{Proportion of Market Share Captured (Exponential Decay)}\\						
	\hline						
	1&	0.12587&	0.12027&	0.14147&	0.16591&	0.19036&	0.21481\\
	2&	0.21392&	0.22854&	0.25671&	0.29059&	0.32447&	0.35836\\
	3&	0.28565&	0.30169&	0.34026&	0.38755&	0.43484&	0.48212\\
	4&	0.35040&	0.37390&	0.41560&	0.47292&	0.53026&	0.58759\\
	5&	0.40132&	0.42225&	0.46572&	0.52422&	0.58360&	0.64297\\
	10&	0.59167&	0.60648&	0.63484&	0.66835&	0.71599&	0.78021\\
	15&	0.69176&	0.70212&	0.71876&	0.74397&	0.78142&	0.84016\\
	20&	0.75030&	0.75737&	0.77181&	0.79276&	0.82209&	0.87325\\
	\hline						
	\multicolumn{7}{|c|}{Run Time in Minutes for all 100 Runs (Exponential Decay)}\\						
	\hline						
	1&	3.21&	9.44&	11.90&	12.45&	10.91&	10.64\\
	2&	6.78&	16.94&	24.63&	23.81&	26.05&	25.34\\
	3&	11.62&	29.54&	38.61&	44.68&	43.46&	43.29\\
	4&	19.58&	47.78&	61.85&	72.02&	70.58&	68.97\\
	5&	34.82&	58.75&	90.84&	105.04&	104.37&	100.30\\
	10&	97.91&	217.74&	295.56&	360.82&	393.02&	383.35\\
	15&	254.80&	533.03&	668.99&	844.88&	1005.76&	959.82\\
	20&	374.09&	1034.84&	1227.17&	1548.28&	1913.28&	1850.52\\
			\hline
		\end{tabular}
	\end{center}
\end{table}

\end{document}